\input amstex
\input Amstex-document.sty
\font\boldmathHuge = cmmib10 scaled 1728

\pageno 155

\def\shorttitle{Some New Developments of Realization of Surfaces into $R^3$}
\def\shortauthor{Jiaxing Hong}

\topmatter
\title\nofrills{\boldHuge Some New Developments of Realization of Surfaces into {\boldmathHuge R}$\,^{\text{\boldLarge 3}}$}
\endtitle

\author \Large Jiaxing Hong* \endauthor

\thanks *Institute of Mathematics, Fudan University, Shanghai 200433,
China. E-mail: \linebreak jxhong\@fudan.ac.cn
\endthanks

\abstract\nofrills \centerline{\boldnormal Abstract}

\vskip 4.5mm

{\ninepoint This paper intends to give a brief survey of the developments on realization of surfaces into $ R^3$
in the last decade. As far as  the local isometric embedding is concerned, some results related to the
Schlaffli-Yau conjecture are reviewed. As for the realization of surfaces in the large, some developments on Weyl
problem for positive curvature and an existence result for realization of complete negatively curved surfaces into
$ R^3$, closely related to Hilbert-Efimov theorem,  are  mentioned. Besides, a few results for two kind of
boundary value problems for realization of positive disks into $R^3$ are introduced.

\vskip 4.5mm


}
\endabstract
\endtopmatter

\document

\baselineskip 4.5mm \parindent 8mm

Given a smooth n-dimensional Riemannian manifold $(M^n$,g),
can we find a map
$$\phi :M^n\longrightarrow R^p\text{ such that $\phi^{*}h=g$}$$
where $h$ is the standard metric in $R^p$?  This is a long standing problem in Differential Geometry. The map
$\phi$ is called isometric embedding or isometric immersion if $ \phi$ is embedding or immersion. There are
several very nice surveys. For example, for known results before 1970, particularly obtained by Russian
mathematicians, see [GR], [G] and for results of $n=2,$ see [Y3] [Y4]. Whereas, what development has been made
during the passed decades? In higher dimensional cases, the most important one is to improve the Nash's theorem so
that $( M^n,g)$ has isometric embedding in a Euclidean space of much lower dimension than that given by Nash and
meanwhile, to use the contraction mapping principle in place of the theorem of the complicated hard implicit
functions, see [H] or [GUN]. In contrast to higher dimensional cases it seems that problems of two dimensional
Riemannian manifolds embedded into $R^3$ has attracted much more attention of mathematicians, particularly in the
last decade. The present paper is devoted to a survey on the developments of isometric embedding (or immersion) of
surfaces into $R^3$ in the last decade. Of course the results mentioned in this survey are by no means exhaustive
and depend a lot on the author's taste.
\par Smoothly and isometrically immersing a surface $(M,g)$ into $
R^3$ is
equivalent to finding three smooth $(C^s,s\ge 1)$ functions $x_{\alpha}
:M\mapsto R^1$, $\alpha =1,2,3$
such that
$$g=dx^2_1+dx^2_2+dx^2_3 .\tag 1$$
Although as far as the formulation of (1) be concerned, $C^1$ regularity
is enough and [K] gives a very nice result in this category. In order
to see the role of the curvature of surfaces in the problem
considered here we prefer to assume $s\ge 2$ throughout the present
paper.
In the sequel, sometimes we use $x,y,z$ to denote $x_1,x_2,x_3$.
In a local coordinates near a point $p\in M$, the metric $g$ is of the
form $g=g_{ij}du^idu^j$. Then (1) can be written as follows
$$\frac {\partial x_{\alpha}}{\partial u^i}\frac {\partial x_{\alpha}}{
\partial u^j}=g_{ij}\text{  }i,j=1,2, \tag 2$$
(2) is a system composed of three differential equations of first
order and hence, this is a determine system.
We say that $\vec {r}$ $=(x_1,x_2,x_3)$ is a local smooth $(C^s)$ isometric embedding
in $R^3$ of the
given surface  if $\vec {r}$ $=(x_1,x_2,x_3)$ is a
smooth $(C^s)$ solution to (2) in a neighbourhood of the point $p$ and that $
\vec {r}$
is a global smooth $(C^s)$ isometric embedding (immersion) in $R^3$ if
$\vec {r}$ $=(x_1,x_2,x_3)$ is a smooth $(C^s)$ solution to (1) on $
M$ and, meanwhile,
is an embedding
(immersion) into $R^3$.
This survey consists of three parts. The first and the second
part include some recent developments in local isometric embedding
and
global isometric embedding respectively. The third part contains some
developments on boundary value problems for realization of positive
disks into $R^3.$
\par {\bf Local isometric embedding}. With the aid of the Cauchy-Kowalevsky
theorem Cartan and Janet proved that any n-dimensional analytic
metric  always admits a local analytic isometric embedding in
$R^{s_n}$ with $s_n=n(n+1)/2$. In the smooth category Gromov in
[GR] proved that any n-dimensional $C^{\infty}$ metric always
admits a local smooth isometric embedding in $R^{s_n+n}$. As
$n=2,$ $s_n=5$ and from (2)  the present result looks far away
from the optimal in the smooth category. On the other hand, [P]
proves that any smooth surface always has a local smooth
isometric embedding in $R^4$. In [Y2, No.22] and also in [Y1, No.
54] Yau posed to prove that any smooth surface always has a local
smooth isometric embedding in $R^3$. In this direction, it is Lin
who first made important breakthrough and his results in [LC1]
and [LC2] state
\proclaim{Theorem 1}{\rm (C. S. Lin)} (1) Any
$C^s,$ $s>10$ nonnegatively curved metric always admits a local
$C^{s-6}$ isometric embedding in $R^3.$
\par (2) If $g$ is a $C^s,$ $s>6$ metric and if its curvature $K$
satisfies
$$K(p)=0\text{ and }dK(p)\ne 0$$
then it admits a $C^{s-3}$ isometric embedding in $R^3$ near $p$. \endproclaim
By means of Lin's technique the problem for local isometric embedding
related to nonpositive curvature metric are also solvable in the following
cases.
\par (3) In [IW],  $K$ nonpositive in a neighbourhood of a point $
p$ and $d^2K(p)\ne 0$,
\par (4) In [HO1], $K=h^{2q}K_1$ where $K_1$, $h$ are smooth functions,
$K_1(p)<0$, $dh(p)\ne 0$ and $q$ is an integer.
\par In what follows let us simply explain what the technique ones
use while  attacking the problem of local isometric embedding.
Suppose that $\vec {r}$ $=(x,y,z)$ is a smooth solution to (2)  in a neighbourhood of the
point in question previously. By the Gauss equations we have, in a
local coordinate system,
$$\vec {r}_{ij}=\Gamma^k_{ij}\vec {r}_k+\Omega_{ij}\vec {n}\text{  or }
\nabla_{ij}\vec {r}=\Omega_{ij}\vec {n},i,j=1,2\tag 3$$
where subscripts i,j  and $\nabla_{ij}$ denote Euclidean and
covariant derivatives respectively and $\Omega_{ij}$ the coefficients of the
second fundamental form, $\Gamma^k_{ij}$ the Christoffel symbols with respect to
the metric and $\vec {n}$ the unit normal to $\vec {r}$. For each unit constant vector,
for instance, the unit vector $\vec {k}$ of the z axis, taking the scale product
of $\vec {k}$ and (3) and using the Gauss equations one can get
$$\det (\nabla_{ij}z)=K\det (g_{ij})(\vec {n},\vec {k})^2 .\tag 4$$
Notice that
$$\align (\vec {n},\vec {k})^2&=1-\left(\frac {(\vec {r}_1\times\vec {r}_2
)\times\vec {k}}{|\vec {r}_1\times\vec {r}_2|}\right)^2\\
&=1-g^{ij}z_iz_j=1-|\nabla z|^2.\endalign$$

Inserting the last expression into (4) we deduce the Darboux
equation
$$F(z)=\det (\nabla_{ij}z)-K\det (g_{ij})(1-|\nabla z|^2)=0.\tag 5$$
Obviously each component of $\vec {r}$ does satisfy the Darboux equation.
Conversely,  for each smooth solution $z$ to (5)
satisfying $|\nabla z|<1,$
$\tilde {g}=g-dz^2$ is a smooth flat metric. Therefore in simply connected
domain $\Omega$ we can always find a smooth mapping  $(x,y$): $\Omega$ $
\mapsto R^2$ such
that $dx^2+dy^2=g-dz^2.$ So the realization of a given metric into $
R^3$
is equivalent to finding a smooth (or $C^s$) solution to the Darboux
equation with a subsidiary condition $|\nabla z|<1.$ If $K$ is positive or
negative at the point considered, then (5)
is elliptic or hyperbolic Monge-Ampere equation and the local
solvability  is well known for both of them. But if $K$ vanishes at this
point, the situation is very complicated and so far there has been no
standard way to deal with such kind of Monge-Ampere equation.
Indeed, its linearized operator
$$L_z\xi =\lim_{t\longrightarrow 0}\frac {F(z+t\xi )-F(z)}t=F^{ij}
\nabla_{ij}\xi +2K(\nabla z,\nabla\xi )\tag 6$$
where $F^{ij}=\partial\det (\nabla_{lk}z)/\partial\nabla_{ij}z$. It is easy to see that the type of this
linear differential operator completely depends on $K$. When $K$
vanishes at the point considered, (6) may be degenerate elliptic,
hyperbolic or mixed type and its local solvability is not clear.
Using a regularized operator instead of (6) and the Nash-Moser procedure
Lin succeeded in proving Theorem 1.
But it is still not clear whether there is obstruction for the
local isometric embedding in smooth category even for the
nonpositive curvature metric. Some results [EG] on linear degenerate
hyperbolic operators of second order which have no local solvability
should be noticed. Anyway, to author's knowledge the problem for local
isometric embedding of surfaces into $R^3$ is still open !!

{\bf Global isometric embedding.}  The first result on global isometric embedding of complete surfaces in $R^3$ is
due to Weyl and Lewy for analytic metric and to Nirenberg and Pogorelov for smooth metric. \proclaim{Theorem }
{\rm (Weyl-Lewy Nirenberg-Pogorelov)} Any analytic (smooth ) positive curvature metric defined on $S^2$ always
admits an analytic ( a smooth ) isometric embedding in $R^3$. \endproclaim
\par For noncompact case, for example,  a complete smooth positive
curvature metric defined on $R^2$, this problem was solved by two
Russian mathematicians. Olovjanisnikov first found the weak
isometric embedding based on the Aleksandrov's theory on convex
surfaces and Pogorelov proved the weak solution smooth if the
metric smooth. \proclaim{Theorem} {\rm (Olovjanisnikov-Pogorelov)}
Any smooth complete positive curvature metric defined on $R^2$
admits a smooth isometric embedding in $R^3$.
\endproclaim

\par The next natural development is to consider the realization in $
R^3$ of nonnegatively curved surfaces. Recently [GL] and [HZ]
independently obtained the following result \proclaim{Theorem 2 }
{\rm (Guan-Li, Hong-Zuily)} Any $C^4$nonnegative curvature metric
defined on $S^2$ always admits a $C^{1.1}$ isometric embedding in
$R^3$.\endproclaim Olovjanisnikov-Pogorelov's result on complete
positively curved plane is also extended to the nonngegatively
curved case, see [HO3]. \proclaim{Theorem 3} Any complete
$C^4$nonnegative curvature metric defined on $R^2$ always admits
a $C^{1.1}$ isometric embedding in $ R^3.$ Moreover it is smooth
where the metric is smooth and the curvature positive.
\endproclaim
In this direction a special case is also obtained in [AM].
\par As far as the regularity of isometric embedding be concerned, ones are
interested in the following question. Can we  improve  the regularity of the isometric embedding obtained in
Theorem 2 and Theorem 3 if the metric is smooth $?$ It is very interesting that [IA] gives a $C^{2.1}$ convex
surface which is not $C^3$ continuous but realizes analytic metric on $S^2$ with positive curvature except one
point. On the other hand, Pogorelov gave a $C^{2,1}$ geodesic disk with nonnegative curvature not even admitting a
$C^2$ local isometric embedding in $R^3$ at the center of this disk. Therefore a natural open question is : Does
there exist a $C^{2.\alpha}$ $(0<\alpha <1)$ isometric embedding in $ R^3$ for any sufficiently smooth (even
analytic) nonnegatively curved sphere or plane?
\par The Hilbert theorem  is one of the most important theorems in
3- Euclidean space. This theorem as well as Efimov's
generalization in [EF1] provide a negative answer for the problem
of realization of complete negatively curved surfaces into $R^3.$
\proclaim{Theorem} {\rm (Hilbert-Efimov)} Any complete surface
with negative constant curvature (with curvature bounded above by
a negative constant) has no $C^2$ isometric immersion in
$R^3$.\endproclaim Another result [EF2] also due to  Efimov
should be mentioned.
\proclaim{Theorem} {\rm (Efimov)} Let $M$ be
a smooth complete negatively curved surface with curvature $K$
subject to
$$\sup_M|K|\text{, }\sup_Mgrad(\frac 1{\sqrt {|K|}})\le C\tag 7$$
for some constant $C$. Then $M$ has no $C^2$ isometric immersion in $
R^3$ \endproclaim
Evidently, Efimov's second result yields a necessary condition for a
complete negatively curved surface to embed isometrically in $R^3,$
$$\sup_M|\nabla\frac 1{\sqrt {|K|}}|=\infty\text{ if }\sup_M|K|<\infty .\tag 8$$
Yau posed the following question [Y1, No.57]. \newline
\centerline{Find a nontrivial sufficient condition for a complete
negatively } \centerline{curved surface to embed isometrically in
$R^3.$} He also pointed out that such a nontrivial condition
might be the rate of decay of the curvature at infinity. Recently
some development in this direction has been made in [HO2]. Let
$M$ be a simply connected noncompact complete surface of negative
curvature $K$. By the Hadamard theorem
$exp_{:}T_p(M)\longrightarrow M$ is a global diffeomorphism for
each point $p\in M$ which induces a global geodesic polar
coordinates $(\rho ,\theta )$ on $M$ centered at $p$.

\proclaim{Theorem 4} Suppose that

(a) for some $\delta>0$, $\rho^{2+\delta}|K|$ is decreasing in $\rho$ outside a compact set and that

(b) $\partial^i_{\theta}\ln |K|$, $i=1,2$ and $\rho
\partial_{\rho}\partial_{\theta}\ln |K|$ bounded on $M$.

\noindent Then $M$ admits a smooth isometric immersion in
$R^3$.\endproclaim
\par\noindent {\bf Remark 1. } If $M\in C^{s.1}(s\ge 4)$ and other assumptions in Theorem 4
are fulfilled, then it admits a $C^{s-1.1}$ isometric immersion in $
R^3$.
\par\noindent {\bf Remark 2. } The assumption (a) implies the rate of decay of the
curvature at the infinity
$$|K|\le\frac A{\rho^{2+\delta}},\text{ for a positive constant }A .\tag 9$$
Such a condition on the decay of the curvature at the infinity
is nearly sharp for the existence since if $\delta =0,$ there might be no
existence.  Consider  a
radius symmetric surface $(R^2,g)$ with
the Gaussian curvature
$$K=-\frac A{1+\rho^2}\text{ for some positive constant }A\tag 10$$
where $\rho$ is the distance function from some point. Evidently
$$\sup_{R^2}|\nabla\frac 1{\sqrt {|K|}}|=\frac 1{\sqrt A}<\infty .$$
Therefore Efimov's second theorem tells us that such a complete negatively
curved surface $(R^2,g)$ has no any $C^2$ isometric immersion in $
R^3$
for arbitrary positive
constant $A$. The arguments of [EF1] and [EF2] are very
genuine but ones expect a more analysis proof for Efimov's results in
[EF1]. Anyway,  a result in [EF3] also by Efimov which is much easier
understood, shows
that the negatively curved surface mentioned above has no $C^2$ isometric
immersion if  $A>3$ in (10).
\par {\bf Boundary value problems for isometric embedding in $\bold R^
3$. } Recently, the study on boundary value problems for isometric
embedding of surface in $R^3$ has attracted much attention of
mathematicians. Various formulations of such questions can be
found in [Y4]. To author's knowledge this field has not been
extensively studied.
 Let $g$ be a smooth positive
curvature metric defined on the closed unit disk $\bar {D}$. Throughout the present paper we always call such
surfaces $(\bar {D},g)$ positive disk. According to the classification in [Y4] there are two kinds of boundary
value problems: one is Dirichlet problem and another is Neumann problem. Let us first consider the Dirichlet
problem. \newline \centerline{Given a smooth positive disk $(\bar {D},g)$ and a complete} \centerline{smooth
surface $\Sigma\subset R^3$, can we find an isometric embedding $ \vec {r}$} \centerline{ of $(\bar {D},g)$ in
$R^3$ such that $\vec {r}(\partial D)\subset\Sigma$?} This problem is also raised in [PO1]. As a first step, one
can consider a simple case.
\newline \centerline{Assume that $\Sigma$ is a plane. Give a complete description of
all }
\centerline{isometric embedding of $(\bar {D},g)$ satisfying $\vec { r}(\partial D)\subset\Sigma$. } Assume
that $\Sigma :\{z=0\}$. Then we are faced with the following boundary value problem.
\newline \centerline{D: To
find an isometric embedding $\vec {r}=(x,y,z)$ of the given}
\centerline{positive disk $(\bar {D},g)$ such that $z(\partial
D)= 0.$}
It is easy to see that there is some obstruction for the
existence of solutions to the above boundary value problem.
Suppose that $\vec {r}=(x,y,z)\in\text{$C^2(\bar {D})$}$ is a
solution to the above boundary value problem. Obviously the
intersection $\vec {r}(\partial D)$ of $\vec { r}$ and the plane
$\{z=0\}$ is a $C^2$ planar convex curve. Denoting its curvature
by $\tilde {k}$ we have
$$\tilde {k}^2=k^2_g+k^2_n\tag 11$$
where $k_g$ and $k_n$ are respectively the geodesic curvature and normal
curvature of $\vec {r}(\partial D)$. $k_n$ is positive everywhere since
$(\bar {D},g)$ has positive curvature. Notice that the total curvature of $
\vec {r}(\partial D)$,
as a planar curve, equals $2\pi$. Hence
$$\int_{\vec {r}(\partial D)}|k_g|ds<2\pi .\tag 12$$
This is a necessary condition in order that the above boundary value
problem $D$ can be solvable. Indeed, this necessary condition is not
sufficient for the solvability. In [HO4] there is a smooth positive disk
satisfying (12) but  not admitting any
$C^2$ solution to the problem $D.$  Furthermore this counter example also
shows that too many changes of
the sign of the geodesic
curvature of the boundary will make the problem $D$ unsolvable.
So we distinguish two cases
$$\text{Case a: }k_g>0\text{ on }\partial D\text{ and }\text{Case b: }
k_g<0\text{ on }\partial D .\tag 13$$
\par It should be pointed out that Pogorelov gave the first solution to the problem $
D$ in [PO1] which states that
\proclaim{Theorem} {\rm (Pogorelov)}
Let $(\bar {D},g)$ be a smooth positive disk. Then the boundary
value problem admits a solution $\vec {r}$ $\in C^{\infty} (D)\cap
C^{0.1}(\bar {D})$ provided that the geodesic curvature of
$\partial D$ with respect to the metric $g$ is nonnegative.
\endproclaim
Pogorelov only obtained a local smooth solution, namely, a solution
smooth inside. One wonder that under what conditions
the problem $D$ always admits a global smooth
solution, namely, a solution smooth up to the boundary. Such global
smooth solutions
are obtained by [DE] for Case a  if there exists a global
$C^2$ subsolution $\psi$ for the Darboux equation (5) in the unit disk $
D$,
vanishing  on
$\partial D$ and with $|\nabla\psi |$ strictly less than 1 on $\bar {
D}$. Recently, [HO5]
removes this technique requirement.
\proclaim{Theorem 5}For Case a the problem $D$ always admits a unique solution
in $C^{\infty}(\bar {D})$ if any one of the following assumptions is satisfied
(1) $K>0$ on $\bar {D}$, (2) $K>0$ in $D$ and $K=0\ne |dK|$ on $\partial
D$. \endproclaim
Obviously, the necessary condition (12) is always satisfied for Case a.
As for Case b it looks rather complicated since there are some
smooth convex surfaces in $R^3$  which are  of no infinitesimal rigidity.
 The presence of such convex surfaces makes  us fail
to prove the existence of the problem
$D$ of Case b by means of the standard method of continuity.
 Thus the solvability of
the problem $D$ for Case b is still open !
\par Let us consider a
spherical crown, in the spherical coordinates,
$$\text{$\Sigma_{\theta_{*}}=\{$}(\sin \theta\cos \phi ,\sin \theta\sin
\phi ,-\cos \theta )|0\le\phi\le 2\pi ,0\le\theta\le\theta_{*}\}$$
and  $\theta =0$ stands for the South pole.  $\Sigma_{\theta_{*}}$
is the isometric embedding of  the metric
$g=d\theta^2+\sin^2\theta d\phi^2,0\le\theta\le\theta_{ *}.$ If
$\theta_{*}>\frac {\pi}2$, then $\Sigma_{\theta_{*}}$ contains
the below hemisphere and the geodesic curvature of its boundary is
negative. We have in [HO4] \proclaim{Theorem 6} There is a
countable set $\Lambda =\{\theta_1,\theta_2,...\theta_n,...\}$
$\subset (\pi /2, \pi )$ with a limit point $\pi /2$ such that
$\Sigma_{\theta_{*}}$ is not infinitesimally rigid if
$\theta_{*}\in \Lambda$. \endproclaim In what follows we proceed
to discuss the Neumann problem for realization of surfaces into
$R^3$. The formulation is as follows. Give  a smooth positive disk
$(\bar {D},g)$ and a positive function $ h\in C^{\infty}(\partial
D)$, The Neumann problem (later, called the problem $N$) is as
follows.
$$N:\text{ }\text{Find a surface }\vec {r}:\bar {D}\mapsto R^3\text{ such that }
d\vec {r}^2=g$$
$$\text{  with the prescribed mean curvature}\text{ }h\text{ on }\vec {
r}(\partial D). \tag 14$$ Let us first introduce an invariant
related to umbilical points of surfaces in $R^3.$ Suppose that
the given metric is of the form
$$\text{$g=Edx^2+2Fdxdy+Gdy^2$,  $(x,y)\in\bar {D}$.}\tag 15$$
Let $\vec {r}$ be a smooth isometric immersion
of $(\bar {D},g)$ with the second fundamental form
$$\text{$II=Ldx^2+2Mdxdy+Ndy^2$}\text{ for }(x,y)\in\bar {D}.\tag 16$$
\proclaim{Definition} If $\vec {r}$ is of no umbilical points on $
\partial D,$ with
$$\sigma =(EM-FL)+\sqrt {-1}(GL-EN)$$
the winding number of $\sigma$ on $\partial D$ is called the index of the umbilical
points of the surface $\vec {r}$ and denoted by $Index(\vec {r})$. \endproclaim
Obviously, this definition makes sense since $p$ is an umbilical point if
and only if $\sigma (p)=0$. Moreover, the definition of
the index of the umbilical points is coordinate-free and hence, an
invariant of describing umbilical points of surfaces. Such invariance
comes from that of a  differential form. Indeed, assume
that in some orthonomal frame, the induced metric and the second
fundamental form of a given surface $\vec {r}$ in $R^3$ are of the form
$$g=\omega^2_1+\omega^2_2\text{ and }II=h_{ij}\omega_i\omega_j$$
respectively. As is well known,
$$\left[(h_{11}-h_{22})+2ih_{12}\right](\omega^2_1+\omega^2_2)$$
is an invariant differential form. So the index of umbilical points is
nothing else but
$$Index(\vec {r})=Index\left\{2h_{12}+i(h_{11}-h_{22})\right\}$$
if no umbilical point on $\partial D$ occurs.
\par The boundary value problem for realization of positive disks into $
R^3$ seems
to have some obstruction. Indeed, even if no imposing any restriction
on the boundary the problem of realization of positive disks into $
R^3$ is
not always solvable.
 For details, refer to
Gromov's counter example [GR] which contains an analytic positive
disk not admitting any $C^2$ isometric immersion in $R^3.$
Therefore for the Neumann boundary value problem
the following hypothesis  is natural. Assume that
$$(\bar {D},g)\text{ admits a $C^2$ isometric immersion $\vec {r}_
0$ in }R^3 .\tag 17$$ We have in [HO6] \proclaim{Theorem 7}  If
$(\bar {D},g)$ is a smooth positive disk satisfying (17) then for
any nonnegative integer $n$ and arbitrary $(n+1)$ distinct points
$p_0\in\partial D,$  $p_1,...,$ $p_n\in D,$ the problem $N$
admits two and only two solutions $\vec {r}$ in $C^{\infty}(\bar
{D},R^3)$ with prescribed mean curvature $ h$ on $\partial D$ and
moreover,
$$\text{one principal direction at }p_0\text{ is tangent to }\partial
D,$$
$$Index(\vec {r})=n\text{ and }H(p_k)=H_0(p_k),k=1,...,n$$
where $H$  and $H_0$ are respectively the mean curvature of $\vec {
r}$ and $\vec {r}_0$
provided that
$$\frac h{\sqrt K}-1>4\max_{\partial D}\left[\frac {H_0}{\sqrt K}-
1\right]\text{ on }\partial D .\tag 18$$
\endproclaim
It is worth pointing out two extreme cases. The first one involves
the existence.
Suppose that the given positive disk $(\bar {D},g)$ is of
positive constant curvature. Then it is easy to see that this positive
disk admits a priori smooth isometric embedding $\vec {r}_0$ in $R^
3$ which is a simply
connected region of the sphere. Under the present circumstance $\vec {
r}_0$ is
totally umbilical and hence, the right hand side of (18) vanishes.
Therefore if $(\bar {D},g)$ is of constant curvature and
$\sqrt K$ $<h\in C^{\infty}(\partial D)$, then the problem $N$ is always solvable for each
nonnegative integer $n$ and arbitrary $(n+1)$ distinct points $p_0
\in\partial D,p_1,,..,p_n\in D.$
\par The second extreme case  involves the nonexistence. If the given positive disk is
radius symmetric, i.e., $g=dr^2+G^2(r)d\theta^2$ $0\le r\le 1$ where $
G\in C^{\infty}([0,1])$
and $G(0)=0,$ $G^{\prime}(0)=1,$ $G>0$ as $r>0.$ Then if $G_r>-1,$ $
(\bar {D},g)$ has such
a priori
smooth isometric embedding in $R^3,$
$$\vec {r}_0:x=G(r)\cos \theta ,y=G(r)\sin \theta ,z=-\int_r^1\sqrt {
1-G^2_r}dr .\tag 19$$ With its mean curvature $H_0=H_0(r)$, if
$H_0(1)>\sqrt {K(1)}$ , then [HO6] proves that for arbitrary
$h\in C^{\infty}(\partial D)$ satisfying $\sqrt {K(1)}\le
h<H_0(1)$ the problem $ N$ has no any $C^2$ solution. Of course,
if $h>4H_0(1)-3\sqrt {K(1)}$, by Theorem 7 the problem $ N$
always admits two and only two smooth solutions for any
nonnegative integer $n$ and arbitrary $n+1$ points
$p_0\in\partial D,$ $p_1,...,p_j\in D$.

\specialhead \noindent \boldLARGE References \endspecialhead

\roster
\item "[AL]" Aleksandrov, A.D., Intrinsic geometry of
convex surfaces, Moscow, 1948, German transl., Academie Verlag
Berlin, 1955.
\item "[AM]" Amano, K., Global isometric embedding of a Riemannian
2-manifold with nonnegative curvature into a Euclidean 3-space,
J.Diff.Geometry, {\bf 35}(1991), 49--83.
\item "[DE]" Delanoe, Ph., Relations globalement regulieres de
disques strictement convexes dans les espaces d'Euclide et de
Minkowski par la methode de Weingarten, Ann.Sci.Ec
Norm.Super.,IV. {\bf 21}(1988), 637--652.
\item "[EF1]" Efimov, N.V., Generalization of singularities on surfaces of
negative curvature, Mat. Sb., 64(1964), 286--320.
\item "[EF2]" Efimov, N.V., A Criterion of homemorphism for some
mappings and its applications to the theory of surfaces, Mat. Sb.,
76(1968), 499--512.
\item "[EF3]" Efimov, N.V., Surfaces with slowly varying negative curvature,
Russ. Math. Survey, 21(1966), 1--55.
\item "[EG]" Egorov, Y., Sur un exemple d'equation lineaire hyperbolique du
second ordre n'ayant pas de solution, Journees  Equations aux
Derivees Partielles, Saint-Jean-de-Monts, 1992.
\item "[G]" Gromov. M., Partial differential relations, Springer-Verlage,
Berlin Heidelberg, 1986.
\item "[GR]" Gromov, M. and Rohklin, V.A., Embeddings and
Immersions in Riemannian Geometry, Russian Mathematical Surveys,
25(1970), 1--57.
\item "[GUN]" Gunther, M., Isometric embeddings of Riemannian manifolds,
Proceeding of the International Congress of Mathematicians, Kyoto,
Japan, 1990.
\item "[GL]" Guan, P. and Li, Y., The weyl problem with nonnegative
Gauss curvature, J. Diff. Geometry, 39(1994), 331--342.
\item "[H]" Hormander, L., On the Nash-Moser implicit function theorem,
Ann. Acad.Sci of Fenn., 10, 255--259.
\item "[HI]" Hilbert, D., Uber Plachen von konstanter Gausscher krummung,
Trans. Amer. Math. Soc., 2(1901),87--99.
\item "[HO1]" Hong, J.X., Cauchy problem for degenerate hyperbolic
Monge-Ampere equations, J.Partial Diff. Equations, {\bf 4}(1991),
1--18.
\item "[HO2]" Hong, J.X., Realization in $R^3$ of complete Riemannian
manifolds with negative curvature, Communications  Anal. Geom.,
1(1993), 487--514.
\item "[HO3]" Hong, J.X., Isometric embedding in $R^3$ of complete
noncompact nonnegatively curved surfaces, Manuscripta Math.,
94(1997), 271--286.
\item "[HO4]" Hong, J.X., Recent developments of realization of surfaces
in $R^3,$AMS/IP, Studies in Advanced Mathematics, Vol.20, (2001),
47--62.
\item "[HO5]"  Hong, J.X., Darboux equations and isometric embedding of
Riemannian manifolds with nonnegative curvature in $R^3$, Chin.
Ann. of math., 2(1999), 123--136.
\item "[HO6]"  Hong, J.X., Positive disks with prescribed mean
curvature on the boundary, Asian J.  Math., Vol.5 (2001), 473-492.
\item "[HZ]" Hong, J.X. and Zuily, C., Isometric embedding
of the 2-sphere with non negative curvature in $R^3$, Math.Z.,
219(1995), 323--334.
\item "[IA]" Iaia, J.A., Isometric embedding of surfaces with
nonnegative curvature in $R^3$, Duke Math. J., 67(1992) 423--459.
\item "[IW]" Iwasaki, N., The stronger hyperbolic equation and its
applications, Proc. Sym. Pure Math., 45(1986), 525--528.
\item "[K]" Kuiper, N.H., On $C^1$-isometric embeddings, I, II, Indag. Math.,
(1955), 545--556, 683--689.
\item "[LC1]" Lin, C.S., The local isometric embedding in $R^3$ of
2-dimensional Riemannian manifolds with nonnegative curvature, J.
Diff. Geometry, 21(1985), 213--230.
\item "[LC2]" Lin, C.S., The local isometric embedding in $R^3$ of two
dimensional Riemannian manifolds with Gaussian curvature changing
sign clearly, Comm. Pure Appl.Math., 39(1986), 307--326.
\item "[NI1]" Nirenberg, L., The Weyl and Minkowski problems in
Differential Geometry in the large, Comm. Pure Appl.Math.,
6(1953), 337--394.
\item "[NI2]" Nirenberg, L., In Nonlinear problems, editor by  R.E.Langer,
University of Wisconsin Press, Madison, 1963, 177--193.
\item "[P]" Poznjak, E., Isometric immersion of 2-dimensional metrics
into Euclidean spaces, Uspechy 28(1970), 47--76.
\item "[PO1]"  Pogorelov, A.V., Extrinsic geometry of convex
surfaces ( transl.  Math. Monogr. Vol 35), Providence, RI,
Am.Math.Soc., 1973.
\item "[PO2]"  Pogorelov, A.V., An example of a two dimensional
Riemannian metric not admitting a local realization in $R^3$,
Dokl. Akad. Nauk. USSR., 198(1971), 42--43.
\item "[Y1]" Yau, S.T., Problem Section, Seminar on Differential
Geometry, Princeton University Press, 1982.
\item "[Y2]" Yau, S.T., Open problems in Geometry, Chern-A great
Geometer on the Twentieth Century, International Press, 1992.
\item "[Y3]" Yau, S.T., Lecture on Differential Geometry, in Berkeley,
1977.
\item "[Y4]" Yau, S.T., Review on Geometry and Analysis, Asian J. of
Math., 4(2000), 235--278.
\endroster

\enddocument